\documentclass[12pt,leqno]{amsart}

\usepackage{amssymb}

\input{diagrams}

\def\sD{\mbox{\sf D}}

\def\D{\sD}

\def\depth{\operatorname{depth}}
\def\dim{\operatorname{dim}}
\def\Ext{\operatorname{Ext}}

\def\fin{\sD^{\operatorname{f}}}

\def\H{\operatorname{H}}

\def\Hom{\operatorname{Hom}}
\def\id{\operatorname{id}}

\def\LTensor{\stackrel{\operatorname{L}}{\otimes}}

\def\RHom{\operatorname{RHom}}

\def\skewtimes{\ltimes\!}

\numberwithin{equation}{part}

\newtheorem{Lemma}{Lemma}[section]
\newtheorem{Theorem}[Lemma]{Theorem}

\newtheorem{Corollary}[Lemma]{Corollary}

\theoremstyle{definition}
\newtheorem{Definition}[Lemma]{Definition}

\newtheorem{Remark}[Lemma]{Remark}

\begin{document}

\title[Recognizing dualizing complexes]
{Recognizing dualizing complexes}

\author{Peter J\o rgensen}
\address{Danish National Library of Science and Medicine, N\o rre
All\'e 49, 2200 K\o \-ben\-havn N, DK--Denmark}
\email{pej@dnlb.dk, www.geocities.com/popjoerg}



\begin{abstract} 
Let $A$ be a noetherian local commutative ring and let $M$ be a
suitable complex of $A$-modules.  This paper proves that $M$ is a
dualizing complex for $A$ if and only if the trivial extension $A
\skewtimes M$ is a Gorenstein Differential Graded Algebra.

As a corollary follows that $A$ has a dualizing complex if and only if
it is a quotient of a Gorenstein local Differential Graded Algebra.
\end{abstract}

\maketitle

\setcounter{section}{0}

Let $A$ be a noetherian local commutative ring and let $M$ be a
complex of $A$-modules with homology $\H\!M$ non-zero and finitely
generated and $\H_i\!M = 0$ for $i < 0$.  Theorem \ref{thm:main1}
shows that $M$ is a dualizing complex for $A$ if and only if the
trivial extension $A \skewtimes M$ is a Gorenstein Differential Graded
Algebra (DGA).  Phrased as a slogan:  DGAs can be used to recognize
dualizing complexes.

In corollary \ref{cor:main2} this is used to show that $A$ has a
dualizing complex if and only if it is a quotient of a Gorenstein
local DGA.

The notion of Gorenstein DGA I shall use is the one from \cite{FJiia};
it is recalled in definition \ref{def:Gorenstein_DGAs}. But note that
for the DGAs in theorem \ref{thm:main1} and corollary \ref{cor:main2},
the condition of being Gorenstein can also be expressed by the familar
equation $\dim_{\ell} \Ext_R(\ell,R) = 1$, see remark
\ref{rmk:Gorenstein_conditions}.  DGAs satisfying this equation were
considered at length in \cite{AFAmerJ}. 

A brief introduction to the theory of DGAs is in \cite{FJiia}.

\setcounter{section}{0}
\section{Definitions}
\label{sec:definitions}

When $A$ is a noetherian commutative ring, $\D(A)$ denotes the derived
category of complexes of $A$-modules, and $\fin(A)$ denotes the full
subcategory of complexes $M$ such that $\H\!M$ is a finitely generated
module over $A$.  The following definition is due to
\cite[def., p.\ 258]{HartsResDual}.

\begin{Definition}
[Dualizing complexes]
\label{def:dualizing_complexes}
Let $A$ be a noetherian local commutative ring.  The complex $D$ in
$\fin(A)$ is called a {\em dualizing complex} for $A$ if the canonical
morphism
\[
  A \longrightarrow \RHom_A(D,D) 
\]
is an isomorphism, and $D$ has finite injective dimension.
\end{Definition}

\medskip
If a dualizing complex exists, then it is unique up to suspension as
follows easily from \cite[thm.\ V.3.1]{HartsResDual}, but existence is
delicate, see \cite[sec.\ 10]{HartsResDual}.

\begin{Definition}
[Trivial extensions]
\label{def:extension}
Let $R$ be a commutative DGA and let $M$ be a Differential Graded
$R$-module (DG $R$-module).  Then $R \oplus M$ is again a DG
$R$-module with differential
\[
  \partial^{R \oplus M} \left( \begin{array}{c} r \\ m \end{array} \right)
  = \left( \begin{array}{c} \partial^R r \\ \partial^M m \end{array} \right),
\]
and the product
\[
  \left( \begin{array}{c} r_1 \\ m_1 \end{array} \right) 
  \cdot \left( \begin{array}{c} r_2 \\ m_2 \end{array} \right)
  = \left( \begin{array}{c} r_1 r_2 \\ r_1 m_2 + m_1 r_2 \end{array} \right)
\]
turns $R \oplus M$ into a DGA called the {\em trivial extension}
of $R$ by $M$, denoted $R \skewtimes M$.
\end{Definition}

\medskip
Observe that there are canonical morphisms of DGAs,
\[
  \begin{array}{ccccc}
    R & \longrightarrow & R \skewtimes M & \longrightarrow & R, \\[.3cm]
    r & \longmapsto
      & \left( \begin{array}{c} r \\ 0 \end{array} \right) \lefteqn{,}
      & 
      & \\[.6cm]
      &
      & \left( \begin{array}{c} r \\ m \end{array} \right) 
      & \longmapsto 
      & r,
  \end{array}
\]
and that the second is a surjection whose kernel is $M$ which can be
viewed as a square zero differential graded ideal in $R \skewtimes M$.

When $R$ is a commutative DGA with $\H_0\!R$ a noetherian ring,
$\D(R)$ denotes the derived category of DG $R$-modules, and $\fin(R)$
denotes the full subcategory of DG modules $M$ such that $\H\!M$ is a
finitely generated module over $\H_0\!R$.  (This is compatible with
the use of the notation $\fin$ given before definition
\ref{def:dualizing_complexes}.)  For commutative DGAs, the definitions
of Gorenstein DGAs and dualizing DG modules from \cite{FJiia} and
\cite{FIJ} simplify as follows.

\begin{Definition}
[Gorenstein DGAs]
\label{def:Gorenstein_DGAs}
Let $R$ be a commutative DGA with $\H_0\!R$ a noetherian ring. 
Then $R$ is cal\-led {\em Gorenstein} if it satisfies:
\begin{enumerate}

  \item  For $M$ in $\fin(R)$, the following biduality morphism is an
         isomorphism,
\[
  M \longrightarrow \RHom_R(\RHom_R(M,R),R).
\]

  \item  The functor $\RHom_R(-,R)$ sends $\fin(R)$ to itself.

\end{enumerate}
\end{Definition}

\begin{Definition}
[Dualizing DG modules]
Let $R$ be a commutative DGA with $\H_0\!R$ a noetherian ring.  The DG
$R$-module $E$ is called a {\em dualizing DG module} for $R$ if it
satisfies:
\begin{enumerate}

  \item  The canonical morphism
\[
  R \longrightarrow \RHom_R(E,E)
\]
is an isomorphism.

  \item  For $M$ in $\fin(R)$ and for $L$ equal to either $R$ or $E$,
         the following evaluation morphism is an isomorphism,
\[
  M \LTensor_R \RHom_R(L,E) \longrightarrow \RHom_R(\RHom_R(M,L),E).
\]

  \item  The functor $\RHom_R(-,E)$ sends $\fin(R)$ to itself.

\end{enumerate}
\end{Definition}

\medskip
It is clear that $R$ is a Gorenstein DGA if and only if it is a
dualizing DG module for itself.

\begin{Definition}
[Local DGAs]
Let $R$ be a DGA.  Then $R$ is called {\em local} if it
satisfies: 
\begin{enumerate}

  \item  $R$ is commutative and concentrated in non-negative
         homological degrees.

  \item  $\H_0\!R$ is a noetherian local ring, and $\H\!R$ is a
         finitely generated module over $\H_0\!R$.

  \item  $R_0$ is a noetherian ring.

\end{enumerate}
The residue class field $\ell$ of $\H_0\!R$ can then be viewed as a
DG $R$-module concentrated in degree $0$, and as such is referred to
as the {\em residue class field} of $R$.
\end{Definition}

\section{Results}
\label{sec:results}

\begin{Lemma}
\label{lem:coinduce}
Let $A$ be a noetherian local commutative ring, let $R$ be a
commutative DGA with $\H_0\!R$ a noetherian ring, and let $R
\longrightarrow A$ be a morphism of DGAs which induces a surjection
$\H_0\!R \longrightarrow A$.

If $R$ is Gorenstein, then
\[
  D = \RHom_R(A,R)
\]
is a dualizing complex for $A$.
\end{Lemma}

\begin{proof}
Observe that the morphism $R \longrightarrow A$ can be used to view
any complex of $A$-modules ${}_{A}M$ as a DG $R$-module ${}_{R}M$. In
other words, any ${}_{A}M$ in $\D(A)$ can be viewed as ${}_{R}M$ in
$\D(R)$.  As $\H_0\!R \longrightarrow A$ is surjective, it is clear
that I have
\begin{equation}
\label{equ:finiteness}
  {}_{A}M \in \fin(A) 
  \; \; \Leftrightarrow \; \; {}_{R}M \in \fin(R).
\end{equation}

First, equation \eqref{equ:finiteness} implies that $A$ viewed
over $R$ is in $\fin(R)$.  Since $R$ is Gorenstein, $D = \RHom_R(A,R)$
is then also in $\fin(R)$, and by equation \eqref{equ:finiteness} this
shows that $D$ is in $\fin(A)$.

Secondly, there are canonical isomorphisms
\begin{align*}
  \RHom_A(D,D)
  & \stackrel{\rm (a)}{=} \RHom_A(\RHom_R(A,R),\RHom_R(A,R)) \\
  & \stackrel{\rm (b)}{\cong} \RHom_R(A \LTensor_A \RHom_R(A,R),R) \\
  & \cong \RHom_R(\RHom_R(A,R),R) \\
  & \stackrel{\rm (c)}{\cong} A,
\end{align*}
where (a) is by the definition of $D$ and (b) is by adjointness, while
(c) is because $R$ is Gorenstein.

Thirdly, let $k$ be the residue class field of $A$. There are
isomorphisms
\begin{align}
  \RHom_A(k,D)
\nonumber
  & \stackrel{\rm (d)}{=} \RHom_A(k,\RHom_R(A,R)) \\
\nonumber
  & \stackrel{\rm (e)}{\cong} \RHom_R(A \LTensor_A k,R) \\
\label{equ:idD}
  & \cong \RHom_R(k,R),
\end{align}
where again (d) is by the definition of $D$ and (e) is by adjointness.
Equation \eqref{equ:finiteness} implies that $k$ viewed over $R$ is in
$\fin(R)$, and since $R$ is Gorenstein, $\RHom_R(k,R)$ is then also in
$\fin(R)$ so has bounded homology. By equation \eqref{equ:idD}, the
same holds for $\RHom_A(k,D)$, and then by
\cite[(A.5.7.4)]{WintherGorDim} the injective dimension $\id_A D$ is
finite.

Altogether, $D$ is a dualizing complex for $A$; cf.\ definition
\ref{def:dualizing_complexes}.
\end{proof}

\begin{Theorem}
\label{thm:main1}
Let $A$ be a noetherian local commutative ring and let $M$ in
$\fin(A)$ have $M \not\cong 0$ and $\H_i\!M = 0$ for $i <
0$. Then
\[
  \mbox{ $M$ is a dualizing complex for $A$ 
         $\Leftrightarrow$ $A \skewtimes M$ is a Gorenstein DGA. }
\]
\end{Theorem}

\begin{proof}
$\Rightarrow$:  Suppose that $M$ is a dualizing complex for $A$.

Clearly, $A \skewtimes M$ can be viewed as
a DGA over $A$ with $\H(A \skewtimes M)$ finitely generated over
$A$. Hence \cite[prop.\ 2.6]{FIJ} says that $A \skewtimes M$ has the
dualizing DG module $\RHom_A(A \skewtimes M,M)$. 

Let 
\[
  M \stackrel{\rho}{\longrightarrow} I
\]
be an injective resolution. Then $\RHom_A(A \skewtimes M,M)$ is
isomorphic to 
\[
  E = \Hom_A(A \skewtimes M,I), 
\]
so $E$ is a dualizing DG module for $A \skewtimes M$.  The $(A
\skewtimes M)$-structure of $E$ comes from the $A \skewtimes M$
appearing in the first variable of the $\Hom$; that is, if $\epsilon$
is a graded element of $E$ and $r_1$ and $r_2$ are graded elements of
$A \skewtimes M$, then
\begin{equation}
\label{equ:R_structure}
  (r_1\epsilon)(r_2) 
  = (-1)^{|r_1|(|\epsilon| + |r_2|)}\epsilon(r_2 r_1).
\end{equation}

Now note that I have
\begin{equation}
\label{equ:columns}
  E = \Hom_A(A \skewtimes M,I) \cong \Hom_A(A,I) \oplus \Hom_A(M,I)
\end{equation}
as complexes of $A$-modules.  This enables me to compute henceforth as
if the elements of $E$ were column vectors $\left( \begin{array}{c}
\alpha \\ \mu \end{array} \right)$ where $A
\stackrel{\alpha}{\longrightarrow} I$ and $M 
\stackrel{\mu}{\longrightarrow} I$ are $A$-linear.  Of course, I must
remember that $E$ is a DG $(A \skewtimes M)$-module via equation
\eqref{equ:R_structure}.

The element $\left( \begin{array}{c} 0 \\ \rho \end{array}
\right)$ in $E$ can be used to define a morphism of
DG $(A \skewtimes M)$-modules by
\begin{equation}
\label{equ:varphi}
  A \skewtimes M \stackrel{\varphi}{\longrightarrow} E, \; \; \;
  \varphi( \left( \begin{array}{c} a \\ m \end{array} \right) )
    = \left( \begin{array}{c} a \\ m \end{array} \right)
      \cdot \left( \begin{array}{c} 0 \\ \rho \end{array} \right).
\end{equation}
In fact, this turns out to be an isomorphism when viewed in $\D(A
\skewtimes M)$. Hence $A \skewtimes M$ is a dualizing DG module for
itself, and so $A \skewtimes M$ is Gorenstein as desired.

To see this, one uses \eqref{equ:R_structure} to get the second $=$ in

\[
  \varphi( \left( \begin{array}{c} a \\ m \end{array} \right) )
  = \left( \begin{array}{c} a \\ m \end{array} \right)
    \cdot \left( \begin{array}{c} 0 \\ \rho \end{array} \right)
  = \left( \begin{array}{c} \chi_{\rho(m)} \\ a \rho \end{array} \right),
\]
where $\chi_i$ denotes the morphism $A \longrightarrow I$ defined by
$\chi_i(a) = ai$, for any $i$ in $I$.  This shows that viewed in
$\D(A)$, the morphism $\varphi$ is the obvious morphism which
identifies $A \oplus M$ with
\begin{align*}
  E 
  & \cong \Hom_A(A,I) \oplus \Hom_A(M,I) \\
  & \cong \RHom_A(A,M) \oplus \RHom_A(M,M) \\
  & \cong M \oplus A.
\end{align*}
So $\varphi$ is an isomorphism when viewed in $\D(A)$.  Hence
$\H\!\varphi$ is bijective, and so $\varphi$ is also an isomorphism
when viewed in $\D(A \skewtimes M)$.

$\Leftarrow$: Suppose that $A \skewtimes M$ is a Gorenstein DGA.

Let me replace $M$ with a quasi-isomorphic complex which consists of
finitely generated modules and satisfies $M_i = 0$ for $i < 0$.
Thereby $A \skewtimes M$ is replaced with a quasi-isomorphic DGA which
is local.  It is clear that the new $A \skewtimes M$ remains
Gorenstein, and that it is enough to show that the new $M$ is a
dualizing complex for $A$.

By the remarks after definition \ref{def:extension}, there is a
canonical morphism $A \skewtimes M \longrightarrow A$ of DGAs, and it
is clear that this induces a surjection $\H_0(A \skewtimes M)
\longrightarrow A$. So since $A \skewtimes M$ is a Gorenstein DGA,
lemma \ref{lem:coinduce} gives that $A$ has a dualizing complex $D$.

Now, \cite[prop.\ 2.6]{FIJ} gives that $A \skewtimes M$ has the
dualizing DG module $\RHom_A(A \skewtimes M,D)$. On the other hand, $A
\skewtimes M$ is Gorenstein, so is a dualizing DG module for itself.
By \cite[thm.\ 3.2]{FIJ} this implies that $\RHom_A(A
\skewtimes M,D)$ and $A \skewtimes M$ are isomorphic up to suspension in
$\D(A \skewtimes M)$. Replacing $D$ with a (positive or negative)
suspension, I can therefore suppose that there is an isomorphism in
$\D(A \skewtimes M)$, 
\[
  A \skewtimes M \stackrel{\cong}{\longrightarrow} 
  \RHom_A(A \skewtimes M,D).
\]
Picking an injective resolution $D \longrightarrow I$ gives that
$\RHom_A(A \skewtimes M,D)$ is isomorphic to
\[
  E = \Hom_A(A \skewtimes M,I)
\]
so all in all there is an isomorphism in $\D(A \skewtimes M)$,
\[
  A \skewtimes M \stackrel{\cong}{\longrightarrow} E.
\]
This can be represented by two quasi-isomorphisms of DG $(A \skewtimes
M)$-modules,
\[
  \begin{diagram}[width=4ex,height=4ex]
    & & X & & \\
    & \ruTo & & \luTo \\
    A \skewtimes M & & & & E \lefteqn{,}
  \end{diagram}
\]

\medskip
\noindent
and as $A \skewtimes M$ is a free DG module over
itself, the left-hand arrow lifts through the right-hand arrow modulo
homotopy, and hence gives a quasi-isomorphism of DG $(A \skewtimes
M)$-modules,
\[
  A \skewtimes M \stackrel{\varphi}{\longrightarrow} E. 
\]

Let me use again the decomposition of $E$ given in equation
\eqref{equ:columns}, and let me suppose
\[
  \varphi( \left( \begin{array}{c} 1 \\ 0 \end{array} \right) )
  = \left( \begin{array}{c} \alpha \\ \mu \end{array} \right).
\]
One can then use \eqref{equ:R_structure} to get the last $=$ in
\begin{align*}
  \varphi( \left( \begin{array}{c} a \\ m \end{array} \right) )
  & = \varphi( \left( \begin{array}{c} a \\ m \end{array} \right)
               \cdot \left( \begin{array}{c} 1 \\ 0 \end{array} \right) ) \\
  & = \left( \begin{array}{c} a \\ m \end{array} \right)
      \cdot \varphi( \left( \begin{array}{c} 1 \\ 0 \end{array} \right) ) \\
  & = \left( \begin{array}{c} \chi_{\alpha(a) + \mu(m)} \\ 
                              a \mu \end{array} \right),
\end{align*}
where again $\chi_i$ denotes the morphism $A \longrightarrow I$
defined by $\chi_i(a) = ai$, for any $i$ in $I$.  This means that if I
view $\varphi$ as a morphism of complexes of $A$-modules,
\[
  A \oplus M 
  \stackrel{\varphi}{\longrightarrow}
  \Hom_A(A,I) \oplus \Hom_A(M,I),
\]
then $\varphi$ is given by the matrix
\[
  \varphi =
  \left(
    \begin{array}{cc}
      \chi_{\alpha(-)} & \chi_{\mu(-)} \\[.15cm]
      \mu \cdot        & 0
    \end{array}
  \right).
\]
The triangular form of the matrix implies that there is a commutative
diagram of complexes of $A$-modules,
\[
  \begin{diagram}[objectstyle=\scriptscriptstyle,labelstyle=\scriptscriptstyle,midshaft,width=4.5ex]
    0 & \rTo & M 
      & \rTo^{\left( \!\!
                     \begin{array}{c}
                       \scriptscriptstyle 0 \\ 
                       \scriptscriptstyle 1
                     \end{array}
                     \!\!
              \right)}
      & A \oplus M 
      & \rTo^{\left( \!\!
                     \begin{array}{cc}
                       \scriptscriptstyle 1 &
                       \scriptscriptstyle 0
                     \end{array}
                     \!\!
              \right)}
      & A & \rTo & 0 \\
    & & \dTo^{\chi_{\mu(-)}} &      
      & \dTo^{\varphi = 
              \left( \!\!\! 
                     \begin{array}{cc}
                       \scriptscriptstyle \chi_{\alpha(-)} \!\!\!\!\!
                       & \scriptscriptstyle \chi_{\mu(-)} \\
                       \scriptscriptstyle \mu \cdot \!\!\!\!\!
                       & \scriptscriptstyle 0
                     \end{array}
                     \!\!\!\!
              \right)}
      & & \dTo_{\mu \cdot} & & \\
    0 & \rTo & \Hom_A(A,I) 
      & \rTo_{\left( \!\!
                     \begin{array}{c}
                       \scriptscriptstyle 1 \\ 
                       \scriptscriptstyle 0
                     \end{array}
                     \!\!
              \right)}
      & \Hom_A(A,I) \oplus \Hom_A(M,I) 
      & \rTo_{\left( \!\!
                     \begin{array}{cc}
                       \scriptscriptstyle 0 &
                       \scriptscriptstyle 1
                     \end{array}
                     \!\!
              \right)}
      & \Hom_A(M,I) & \rTo & 0 \lefteqn{,} \\
  \end{diagram}
\]
which can be written more simply as
\[
  \begin{diagram}[objectstyle=\scriptscriptstyle,labelstyle=\scriptscriptstyle,midshaft]
    0 & \rTo & M 
      & \rTo
      & A \oplus M 
      & \rTo
      & A & \rTo & 0 \\
    & & \dTo & & \dTo^{\varphi} & & \dTo & & \\
    0 & \rTo & I
      & \rTo
      & I \oplus \Hom_A(M,I) 
      & \rTo
      & \Hom_A(M,I) & \rTo & 0 \lefteqn{.} \\
  \end{diagram}
\]

Let me apply to this the functor $\Hom_A(P,-)$ where $P
\longrightarrow k$ is a projective resolution of
$k$, the residue class field of $A$. There results a new commutative
diagram of complexes of $A$-modules, with split exact rows (because
the above diagram has split exact rows), and with the middle vertical
morphism a quasi-isomorphism (because the same holds for $\varphi$ in
the above diagram, and because $\Hom_A(P,-)$ preserves
quasi-isomorphisms).

The new commutative diagram of complexes of $A$-modules induces a
commutative diagram of long exact sequences of homology groups. Since
$P$ is a projective resolution of $k$, the homology groups are certain
$\Ext_A^i(k,-)$'s.  Moreover, the connecting homomorphisms are zero
(because the rows in the diagram of complexes are split exact), and
the vertical morphisms which result from the middle vertical morphism in
the diagram of complexes are isomorphisms (because the middle vertical
morphism in the diagram of complexes is a quasi-isomorphism).

Summing up, this gives for each $i$ a commutative diagram with exact
rows, 
\[
  \begin{diagram}[objectstyle=\scriptscriptstyle,labelstyle=\scriptscriptstyle,midshaft,width=5ex]
    0 & \rTo & \Ext_A^i(k,M) & \rTo & \Ext_A^i(k,A \oplus M)
      & \rTo & \Ext_A^i(k,A) & \rTo & 0 \\
    & & \dTo & & \dTo_{\cong} & & \dTo & & \\
    0 & \rTo & \Ext_A^i(k,I) & \rTo 
      & \Ext_A^i(k,I \oplus \Hom_A(M,I)) & \rTo
      & \Ext_A^i(k,\Hom_A(M,I)) & \rTo & 0 \lefteqn{.} \\
  \end{diagram}
\]

In particular, $\Ext_A^i(k,M) \longrightarrow \Ext_A^i(k,I)$ is
injective for each $i$, so
\[
  \dim_k \Ext_A^i(k,M) \leq \dim_k \Ext_A^i(k,I)
                       =    \dim_k \Ext_A^i(k,D)
\]
holds for each $i$. But I have $\dim_k \Ext_A^i(k,D) = \delta_{ie}$
for some fixed $e$ by \cite[prop.\ V.3.4]{HartsResDual}, where
$\delta_{ie}$ is $1$ for $i = e$ and $0$ otherwise.  So the only
possibilities for $\dim_k \Ext_A^i(k,M)$ are that it is $\delta_{ie}$,
or that it is identically zero.

The latter alternative would give $\depth_A M = \infty$, but this is
impossible by \cite[lem.\ (A.8.9)]{WintherGorDim} because of $M
\not\cong 0$. So the former alternative
\[
  \dim_k \Ext_A^i(k,M) = \delta_{ie}
\]
must hold, and then $M$ is a dualizing complex for $A$ by
\cite[prop.\ V.3.4]{HartsResDual} again.
\end{proof}

\begin{Corollary}
\label{cor:main2}
Let $A$ be a noetherian local commutative ring.  Then $A$ has a
dualizing complex if and only if it is a quotient of a Gorenstein local
DGA.
\end{Corollary}

\begin{proof}
Suppose that $A$ is a quotient of $R$ which is a Gorenstein local DGA,
and let $R \longrightarrow A$ be the quotient morphism.  As $R$ is
concentrated in non-negative degrees, this clearly induces a
surjection $\H_0\!R \longrightarrow A$. So $D = \RHom_R(A,R)$ is a
dualizing complex for $A$ by lemma \ref{lem:coinduce}.

On the other hand, suppose that $A$ has a dualizing complex $D$.  By
replacing $D$ with a high suspension, I can suppose $\H_i\!D = 0$ for
$i < 0$, and then by replacing $D$ with a quasi-isomorphic complex, I
can suppose that $D$ consists of finitely generated modules and
satisfies $D_i = 0$ for $i < 0$.  Then $A \skewtimes D$ is a local
DGA, and as $D$ is a dualizing complex for $A$, theorem
\ref{thm:main1} gives that $A \skewtimes D$ is a Gorenstein DGA.  And
$A \skewtimes D$ has $A$ as a quotient by the remark after definition
\ref{def:extension}.
\end{proof}

\begin{Remark}
\label{rmk:Gorenstein_conditions}
Theorem \ref{thm:main1} and corollary \ref{cor:main2}
consider DGAs which are Gorenstein in the sense of definition
\ref{def:Gorenstein_DGAs}. However, for the DGAs in question, this
condition can be expressed in an alternative, simple way:

In \cite[thm.\ 4.3]{FIJ} was proved that if $R$ is a local DGA with
residue class field $\ell$, then 
\begin{equation}
\label{equ:AF_Gorenstein}
  R \mbox{ is Gorenstein }  \Leftrightarrow \;
  \dim_{\ell} \Ext_R(\ell,R) = 1.
\end{equation}

Corollary \ref{cor:main2} deals with a local DGA to which equation
\eqref{equ:AF_Gorenstein} applies.

Theorem \ref{thm:main1} does not deal directly with a local DGA.
However, let me replace $M$ in theorem \ref{thm:main1} with a
quasi-isomorphic complex which consists of finitely generated modules
and satisfies $M_i = 0$ for $i < 0$.  Thereby $A \skewtimes M$ is
replaced with a quasi-isomorphic DGA which {\em is} local and to which
equation \eqref{equ:AF_Gorenstein} applies.  And it is clear that the
old and the new $A \skewtimes M$ are Gorenstein simultaneously.
\end{Remark}

\noindent
{\bf Acknowledgements. }  Precursors to the above results have been
mentioned in conversations I have had with Amnon Yekutieli and Anders
Fran\-kild. 

The principle of re\-cog\-ni\-zing dualizing complexes by means of
trivial extensions is well known in the Cohen-Macaulay case, see
\cite[prop.\ 4.2]{FoxbyGorMod} and \cite[thm.\ (7)]{ReitenThm}.

By \cite[thm.\ 1.2]{Kawasaki}, a stronger result holds than corollary
\ref{cor:main2}, namely, $A$ has a dualizing complex if and only if it
is a quotient of a Gorenstein noetherian local commutative ring $R$.

A different way of using DGAs to recognize dualizing complexes is in
\cite{HinichDual}.

The diagrams were typeset with Paul Taylor's {\tt diagrams.tex}.

\end{document}